\documentclass{amsart}
\usepackage{cases}
\usepackage{amssymb,amsmath}
\usepackage{pifont}
\usepackage{amsfonts}
\usepackage{txfonts}
\usepackage{mathrsfs}
\usepackage{bbm}
\usepackage{enumerate}
\usepackage{color}

\theoremstyle{plain}
\newtheorem{theorem}{Theorem}[section]
\newtheorem{remark}[theorem]{Remark}
\newtheorem{lemma}[theorem]{Lemma}
\newtheorem{proposition}[theorem]{Proposition}

\numberwithin{equation}{section}

\begin{document}

\title[Non-cutoff Boltzmann equation]
{Shubin regularity for the radially symmetric spatially homogeneous
Boltzmann equation with Debye-Yukawa potential}

\author[L\'eo Glangetas,Hao-Guang Li]
{L\'eo Glangetas,Hao-Guang Li}

\address{L\'eo Glangetas,
\newline\indent
Universit\'e de Rouen, CNRS UMR 6085, Math\'ematiques
\newline\indent
76801 Saint-Etienne du Rouvray, France}
\email{leo.glangetas@univ-rouen.fr}

\address{Hao-Guang Li,
\newline\indent
School of Mathematics and statistics,
\newline\indent
South-central university for nationalities 430074,Wuhan, P. R. China}
\email{lihaoguang@mail.scuec.edu.cn}

\date{\today}



\begin{abstract}
In this work, we study the Cauchy problem for the radially symmetric spatially homogeneous
Boltzmann equation with Debye-Yukawa potential. We prove that this Cauchy problem enjoys
the same smoothing effect as the Cauchy problem defined by the evolution equation associated to a fractional
logarithmic harmonic oscillator.  To be specific, we can prove the solution of the Cauchy problem belongs to Shubin spaces.
\end{abstract}

\maketitle
\tableofcontents


\section{Introduction}\label{S1}
In this work, we consider the spatially homogeneous Boltzmann equation
\begin{equation}\label{eq1.10}
\frac{\partial f}{\partial t} = Q(f,f)
\end{equation}
where $f=f(t,v)$ is the density distribution function depending only on two varia\-bles $t\geq0$~and $v\in\mathbb{R}^{3}$. The Boltzmann bilinear collision operator is given by
\begin{equation*}
Q(g,f)(v)=\int_{\mathbb{R}^{3}}\int_{S^{2}}B(v-v_{\ast},\sigma)\left(g(v_{\ast}^{\prime})f(v^{\prime})-g(v_{\ast})f(v)\right)dv_{\ast}d\sigma,
\end{equation*}
where for $\sigma\in \mathbb{S}^{2}$,~the symbols~$v_{\ast}^{\prime}$~and~$v^{\prime}$~are abbreviations for the expressions,
$$
v^{\prime}=\frac{v+v_{\ast}}{2}+\frac{|v-v_{\ast}|}{2}\sigma,\,\,\,\,\, v_{\ast}^{\prime}
=\frac{v+v_{\ast}}{2}-\frac{|v-v_{\ast}|}{2}\sigma,
$$
which are obtained in such a way that collision preserves momentum and kinetic energy,~namely
$$
v_{\ast}^{\prime}+v^{\prime}=v+v_{\ast},\,\,\,\,\, |v_{\ast}^{\prime}|^{2}+|v^{\prime}|^{2}=|v|^{2}+|v_{\ast}|^{2}.
$$
For monatomic gas, the collision cross section $B(v-v_{\ast},\sigma)$ is a non-negative function which~depends only on $|v-v_{\ast}|$ and $\cos\theta$
which is defined through the scalar product in $\mathbb{R}^{3}$ by
$$\cos\theta=\frac{v-v_{\ast}}{|v-v_{\ast}|}\cdot\sigma.$$
Without loss of generality, we may assume that $B(v-v_{\ast},\sigma)$ is supported on the set $\cos\theta\geq0,$ i.e. where $0\leq \theta\leq\frac{\pi}{2}$.
See for example \cite{NYKC1}, \cite{Villani} for more explanations about the support of $\theta$.
For physical models, the collision cross section usually takes the form
\begin{equation*}
B(v-v_{\ast},\sigma)=\Phi(|v-v_{\ast}|)b(\cos\theta),
\end{equation*}
with a kinetic factor
$$\Phi(|v-v_{\ast}|)=|v-v_{\ast}|^{\gamma},\,\,\gamma\in]-3,+\infty[.$$
The molecules are said to be Maxwellian when the parameter $\gamma=0$.

Except for the hard sphere model, the function $b(\cos\theta)$ has a singularity at $\theta=0.$
For instance, in the important model case of the inverse-power potentials,
$$U(\rho)=\frac{1}{\rho^r}, \,\,\text{with}\,\, r>1,$$
with $\rho$ being the distance between two interacting particles
in the physical 3-dimensional space $\mathbb{R}^{3}$,
\begin{equation*}
b(\cos\theta)\sin\theta\sim K\theta^{-1-\frac{2}{r}}, \,\,\text{as}\,\,\theta\rightarrow0^+.
\end{equation*}
The notation $a\sim b$ means that there exist positive constants $C_2>C_1>0$, such that
$$C_1 \, a\leq b\leq C_2 \, a.$$
Notice that the Boltzmann collision operator is not well defined for the case when $r=1$ corresponding to the Coulomb potential.

If the inter-molecule potential satisfies the Debye-Yukawa type potentials where the potential function is given by
$$U(\rho)=\frac{1}{\rho\, e^{\rho^s}},\,\,\text{with}\,\,s>0,$$
the collision cross section has a singularity in the following form
\begin{equation}\label{b}
b(\cos\theta)\sim \theta^{-2}(\log\theta^{-1})^{\frac{2}{s}-1},\,\,\text{when}\,\,\theta\rightarrow 0^+,  \,\,\text{with}\,\,s>0.
\end{equation}
This explicit formula was first appeared in the Appendix in \cite{YSCT}.
In some sense, the Debye-Yukawa type potentials is a model between the Coulomb potential corresponding to $s=0$ and the inverse-power potential:
{\color{black}
This behavior can be computed from the equations
(conservation of energy and angular
momentum respectively)
\begin{align*}
&\frac12 (\dot{\rho}^2 + \rho^2 \dot{\varphi}^2) + U(\rho) = \frac12 V^2 + U(\sigma),
\\
&\rho^2 \dot{\varphi} = p(V,\tilde{\theta}) V^2
\end{align*}
where $\rho$ and $\varphi$ are the radial and angular coordinates in the plane of motion
and $p(V, \tilde{\theta})$ is the impact parameter which defines the collision cross section
\begin{align*}
B(|z|,\tilde{\theta}) = |z| \frac{p}{2 \sin\tilde{\theta}} \frac{\partial p}{\partial\tilde{\theta}}
\end{align*}
and $z = v - v_*$ is the relative velocity,
$\theta = \pi - 2 \tilde{\theta}$ is the deviation angle.
}
For further details on the physics background and the derivation of the Boltzmann equation, we refer to the references \cite{Cerci}, \cite{Villani}.

We linearize the Boltzmann equation near the absolute Maxwellian distribution
$$
\mu(v)=(2\pi)^{-\frac 32}e^{-\frac{|v|^{2}}{2}}.
$$
Let $f(t,v)=\mu(v)+\sqrt{\mu}(v)g(t,v)$. Plugging this expression into $\eqref{eq1.10}$, we have
$$
\frac{\partial g}{\partial t}+\mathcal{L}[g]={\bf \Gamma}(g, g)
$$
with
$$
{\bf \Gamma}(g, h)=\frac{1}{\sqrt{\mu}}Q(\sqrt{\mu}g,\sqrt{\mu}h),\,\,
\mathcal{L}(g)=-\frac{1}{\sqrt{\mu}}[Q(\sqrt{\mu}g,\mu)+Q(\mu,\sqrt{\mu}g)].
$$
Then the Cauchy problem \eqref{eq1.10} can be rewrited in the form
\begin{equation}\label{eq-1}
\left\{ \begin{aligned}
         &\partial_t g+\mathcal{L}(g)={\bf \Gamma}(g, g),\,\\
         &g|_{t=0}=g_{0}.
\end{aligned} \right.
\end{equation}
The linear operator $\mathcal{L}$ is nonnegative (\cite{NYKC1,NYKC2,NYKC3}),\,
with the null space
$$
\mathcal{N}=\text{span}\left\{\sqrt{\mu},\,\sqrt{\mu}v_1,\,\sqrt{\mu}v_2,\,
\sqrt{\mu}v_3,\,\sqrt{\mu}|v|^2\right\}.
$$
Denote by $\mathbf{P}$ the orthoprojection from $L^2(\mathbb{R}^{3})$
into $\mathcal{N}$.\,\,Then
$$
(\mathcal{L}g,\,g)=0\Leftrightarrow g=\mathbf{P}g.
$$
In the case of the inverse-power potential with $r>1$,
the regularity of the Boltzmann equation has been studied by numerous papers.  Regarding the Cauchy problem \eqref{eq1.10},
it is well known that the non-cutoff spatially homogeneous Boltzmann equation
enjoys an $\mathscr{S}(\mathbb{R}^3)$-regularizing effect
for the weak solutions to the Cauchy problem \eqref{eq1.10}(see \cite{DW,YSCT}).
For the Gevrey regularity, Ukai showed in \cite{Ukai} that the Cauchy problem for the Boltzmann
equation has a unique local solution in Gevrey classes. Then, Desvillettes, Furioli and Terraneo proved in \cite{DFT}
the propagation of Gevrey regularity for solutions of the Boltzmann equation with Maxwellian mole\-cules. For
mild singularities, Morimoto and Ukai proved in \cite{MU} the Gevrey regularity
of smooth Maxwellian decay solutions to the Cauchy problem of the spatially homogeneous
Boltzmann equation with a modified kinetic factor. See also \cite{TZ}
for {\color{black} the non-modified} case.
On the other hand, Lekrine and Xu proved in \cite{L-X} the property of
Gevrey smoothing effect for the weak solutions to the Cauchy problem associated to the radially symmetric spatially
homogeneous Boltzmann equation with Maxwellian molecules for $r>2$.
This result was then completed by Glangetas
and Najeme who established in \cite{G-N} the analytic smoothing effect in the case when $1<r<2$.
In \cite{BHRV0}, it has led to the hope that the homogenous Boltzmann equation enjoys similar regularity properties as the heat equation with a fractional Laplacian.
Regarding the linearized Cauchy problem \eqref{eq-1}, it has been proved that the solutions for linearized non-cutoff Boltzmann
 equation belongs to the symmetric Gelfand-Shilov spaces $S^{r/2}_{r/2}(\mathbb{R}^3)$
for any positive time, see \cite{ NYKC1}, \cite{HAOLI}.
The Gelfand-Shilov space $S^{\nu}_{\nu}(\mathbb{R}^3)$
with $\nu\geq\frac{1}{2}$ can be identify with
$$
S^{\nu}_{\nu}(\mathbb{R}^{3})=\left\{f\in C^\infty (\mathbb{R}^3);  \exists \tau>0,
\|e^{\tau \, \mathcal{H}^{\frac{1}{2\nu}}}f\|_{L^2}<+\infty\right\}.
$$
where $\mathcal{H}$ is the harmonic oscilator
\begin{equation*}
  \mathcal{H}=-\triangle +\frac{|v|^2}{4}.
\end{equation*}
For the Cauchy problem \eqref{eq-1},  it has been proved in  \cite{NYKC3} and \cite{GLX_2015} that the Cauchy problem for the non-cutoff spatially homogeneous Boltzmann equation with the small initial datum $g_0\in L^2(\mathbb{R}^3)$ has a global solution, which belongs to the Gelfand-Shilov class $S^{r/2}_{r/2}(\mathbb{R}^3)$.


In the present work, we consider the collision kernel in the Maxwellian molecules case and the angular function $b$ satisfying the Debye-Yukawa potential \eqref{b} for some $s>0$.
For convenience, we denote
\begin{equation}\label{beta}
\beta(\theta)=2\pi b(\cos\theta)\sin\theta.
\end{equation}
We study the smoothing effect for the Cauchy problem \eqref{eq-1} associated to the non-cutoff
spatially homogeneous Boltzmann equation with Debye-Yukawa potential \eqref{b}.
The singularity of the collision kernel $b$
endows the linearized Boltzmann operator $\mathcal{L}$ with
the logarithmic regularity property, see Proposition 2.1 in \cite{GL-2015}, the linearized Debye-Yukawa potential Boltzmann operator $\mathcal{L}$ was shown to behave as a fractional logarithmic harmonic oscilator $\left(\log(\mathcal{H}+1)\right)^{\frac{2}{s}}$.
The logarithmic regularity theory was first introduced in \cite{Morimoto} on the hypoellipticity of the infinitely degenerate elliptic operator and was developed in \cite{Morimoto-Xu1},\cite{Morimoto-Xu2} on the logarithmic Sobolev estimates.
Recently, for $0<s<2$, in \cite{YSCT} it has been shown that weak solutions to the Cauchy problem \eqref{eq1.10} with Debye-Yukawa
type interactions enjoy an $H^{\infty}$ smoothing property, i.e. starting with arbitrary initial datum $f_0\geq0$,
$$\int_{\mathbb{R}^3}f_0(v)(1+|v|^2+\log(1+f_0(v)))dv<+\infty,$$
one has $f(t,\cdot)\in H^{\infty}(\mathbb{R}^3)$ for any positive time $t > 0$.
This result was extended by J.-M. Barbaroux, D. Hundertmark, T. Ried, S. Vugalter in \cite{BHRV}.
{\color{black} They showed
a stronger 
regularisation property : 
}
for any $0<s<2$,
and for any $T_0>0$, there exist $\beta, M>0$ such that
$$e^{\beta t (\log\langle D_v\rangle)^{\frac{2}{s}}}f(t,\cdot)\in L^2(\mathbb{R}^d)$$
and
$$\sup_{\lim \eta\in\mathbb{R}^d}e^{\beta t (\log\langle D_v\rangle)^{\frac{2}{s}}}|\hat{f}(t,\eta)|\leq M$$
for all $t\in(0,T_0] $ with $\langle v\rangle=(1+|v|^2)^{1/2}.$

{\color{black} In this paper, we improve the regularisation property (for small initial data)}.
Based upon our recent results 
\cite{NYKC3} and \cite{GLX_2015} of the Gelfand-Shilov smoothing effect for the homogeneous Boltzmann equation with Maxwellian molecules in the case of the inverse-power potential and the result of \cite{GL-2015} for the linear homogeneous Boltzmann equation with Debye-Yukawa potential,
 we 
show that, for small initial data,
the Cauchy problem \eqref{eq-1} for the radially symmetric homogeneous Boltzmann equation enjoys the same smoothing effect as the Cauchy problem defined by the evolution equation associated to a fractional
logarithmic harmonic oscillator.
In order to precise the regularizing effect of the solution for the Cauchy problem \eqref{eq-1},
we introduce the Shubin spaces.
Let $\tau\in\mathbb{R}$,
we denote by $Q^{\tau}(\mathbb{R}^3)$ the spaces introduced by Shubin \cite{Shubin}, Ch. IV, 25.3, with norm
\[\|u\|_{Q^{\tau}(\mathbb{R}^3)} =
\Bigl\|\Bigl( -\Delta+{\textstyle \frac{|v|^2}{4} + 1} \Bigr)^{\frac{\tau}{2}} \, u\Bigr\|_{L^2(\mathbb{R}^3)}=\Bigl\|\Bigl(\mathcal{H}+ 1 \Bigr)^{\frac{\tau}{2}} \, u\Bigr\|_{L^2(\mathbb{R}^3)}.\]

Now we begin to present our results.
\begin{theorem}\label{trick}
Assume that the Maxwellian collision cross-section $b(\,\cdot\,)$ is given in~$\eqref{b}$ with $0<s\leq2$,
then there exists $\varepsilon_0>0$ such that
for any initial datum $g_0\in\,L^2(\mathbb{R}^3)\bigcap\mathcal{N}^{\perp}$ with $\|g_0\|^2_{L^2(\mathbb{R}^3)}\le \varepsilon_0$,
the Cauchy problem \eqref{eq-1}~admits  a solution which
belongs to any Shubin spaces for any $t>0$.
Furthermore, there exist  $c_0>0$, $C>0$ such that, for any $t\ge 0$,
\begin{equation}\label{rate}
\|e^{tc_0\left(\log(\mathcal{H}+1)\right)^{\frac{2}{s}}}g\|_{L^2}
\leq\,
Ce^{-\frac{\lambda_{2,0}t}{4} } \, \|g_0\|_{L^2(\mathbb{R}^3)},
\end{equation}
where
$$
\lambda_{2,0}=\int_{|\theta|\leq\pi/4}\beta( \theta)(1-\sin^4\theta-\cos^4\theta)d\theta>0.
$$
To be more specific,\\
1) in the case $0<s\leq2$:
\begin{equation}\label{rate1}
\forall t>0,\quad
\|g(t)\|_{Q^{2c_0t}}\
\leq
 e^{-\frac{\lambda_{2,0}t}{4}}\|g_{0}\|_{L^2(\mathbb{R}^3)}.
\end{equation}
%
\\
2) in the case $0<s<2$,
there exists a constant $c_s>0$ such that for any $t>0$,
\begin{equation}\label{rate2}
\forall k\geq 0,\quad
\|g(t)\|_{Q^{k}}
\leq
 e^{-\frac{\lambda_{2,0}t}{4}}e^{c_s \, (1/t)^{\frac{s}{2-s}}  \, k^{\frac{2}{2-s}}}\|g_{0}\|_{L^2(\mathbb{R}^3)}.
\end{equation}
\end{theorem}

\begin{remark}
We have proved that, if the initial data $g_0$ is small enough and contained in $L^2(\mathbb{R}^3)$ in Cauchy problem \eqref{eq-1} , then the global
solution for the Cauchy problem \eqref{eq1.10} return to the equilibrium with respect to Shubin space norm.
\end{remark}

\begin{remark}
We think that the regularity properties are optimal, since they are optimal
concerning the linearised Cauchy problem (see \cite{GL-2015}.
\end{remark}

The rest of the paper is arranged as follows. In Section \ref{S2}, we introduce the spectral analysis of the linear Boltzmann operator and in Section \ref{S3}, we establish an upper bounded estimates
of the nonlinear operators with an exponential weighted norm.
The proof of the main Theorem \ref{trick} will be presented in Section \ref{S4}.  In Section \ref{section estimate}, we provide the proof of the technical Lemma \ref{sum}.
In the Appendix \ref{appendix}, we present some indentity properties
of the Shubin spaces used in this paper and the proof of some technical Lemmas.

\section{The spectral analysis of the Boltzmann operators }\label{S2}

\subsection{Diagonalization of linear operators}

We first recall the spectral decomposition of linear Boltzmann operator.
In the cutoff case, that is, when $b(\cos\theta)\sin\theta\in\,L^1([0,\frac{\pi}{2}])$, it was shown in \cite{WU} that
\begin{equation*}
\mathcal{L}(\varphi_{n, l, m})=\lambda_{n,l}\, \varphi_{n, l, m}, \,\,n,l\in\mathbb{N},\,m\in\mathbb{Z}, |m|\leq l.
\end{equation*}
This diagonalization of the linearized Boltzmann operator with Maxwellian molecules holds
as well in the non-cutoff case, (see \cite{Boby,Cerci,Dole,NYKC1,NYKC2}).\,\,Where
\begin{equation*}
\lambda_{n,l}=\int_{|\theta|\le \frac{\pi}{4}}\beta(\theta)\Big(1+\delta_{n, 0}\delta_{l, 0}
-(\sin\theta)^{2n+l}P_{l}(\sin\theta)-(\cos\theta)^{2n+l}P_{l}(\cos\theta)\Big)d\theta,
\end{equation*}
the eigenfunctions are
\begin{equation*}%
\varphi_{n,l,m}(v)=\left(\frac{n!}{\sqrt{2}\Gamma(n+l+3/2)}\right)^{1/2}
\left(\frac{|v|}{\sqrt{2}}\right)^{l}e^{-\frac{|v|^{2}}{4}}
L^{(l+1/2)}_{n}\left(\frac{|v|^{2}}{2}\right)Y^{m}_{l}\left(\frac{v}{|v|}\right),
\end{equation*}
where $\Gamma(\,\cdot\,)$ is the standard Gamma function, for any $x>0$,
$$\Gamma(x)=\int^{+\infty}_0t^{x-1}e^{-x}dx.$$
The $l^{th}$-Legendre polynomial~$P_{l}$ and the Laguerre polynomial $L^{(\alpha)}_{n}$~of order $\alpha$,~degree $n$ (see\,\cite{San})\,read,
\begin{align*}
&P_{l}(x)=\frac{1}{2^ll!}\frac{d^l}{dx^l}(x^2-1)^l,\,\,\text{where}\,|x|\leq1;\\
&L^{(\alpha)}_{n}(x)=\sum^{n}_{r=0}(-1)^{n-r}\frac{\Gamma(\alpha+n+1)}{r!(n-r)!
\Gamma(\alpha+n-r+1)}x^{n-r}.
\end{align*}
For any unit vector $\sigma=(\cos\theta,\sin\theta\cos\phi,\sin\theta\sin\phi)$~with $\theta\in[0,\pi]$~and~$\phi\in[0,2\pi]$,~the orthonormal basis of spherical harmonics~$Y^{m}_{l}(\sigma)$ is
\begin{equation*}
Y^{m}_{l}(\sigma)=N_{l,m}P^{|m|}_{l}(\cos\theta)e^{im\phi},\,\,|m|\leq l,
\end{equation*}
where the normalisation factor is given by
$$
N_{l,m}=\sqrt{\frac{2l+1}{4\pi}\cdot\frac{(l-|m|)!}{(l+|m|)!}}
$$
and $P^{|m|}_{l}$~is the associated Legendre functions of the first kind of order $l$ and degree $|m|$ with
\begin{equation}\label{Plm}
P^{|m|}_{l}(x)= (1-x^2)^\frac{|m|}{2}
\left(\frac{\mathrm{d}}{\mathrm{d}x}\right)^{|m|} P_{l}(x).
\end{equation}
The family $\Big(Y^m_l(\sigma)\Big)_{l\geq0,|m|\leq\,l}$ constitutes an orthonormal basis of the space $L^2(\mathbb{S}^2,\,d\sigma)$ with $d\sigma$ being the surface measure on $\mathbb{S}^2$ (see \cite{Jones}, \cite{JCSlater}).  Noting that $\left\{\varphi_{n,l,m}\right\}$ consist an orthonormal basis of $L^2(\mathbb{R}^3)$ composed of eigenvectors of the harmonic oscillator
(see\cite{Boby}, \cite{NYKC2})
\begin{equation*}
\mathcal{H}(\varphi_{n, l, m})=(2n+l+\frac 32)\, \varphi_{n, l, m}.
\end{equation*}
As a special case, $\left\{\varphi_{n, 0, 0}\right\}$ consist an orthonormal basis of $L^2_{rad}(\mathbb{R}^3)$ in the radially symmetric function space (see \cite{NYKC3}) and
\begin{equation*}
\mathcal{H}(\varphi_{n, 0, 0})=(2n+\frac 32)\, \varphi_{n, 0, 0}.
\end{equation*}
We have that, for suitables functions
$g$,
\begin{equation*}
\mathcal{L}(g)=\sum^{\infty}_{n=0}\sum^{\infty}_{l=0}\sum^{l}_{m=-l}
\lambda_{n,l}\, g_{n,l,m}\, \varphi_{n, l, m},
\end{equation*}
where $g_{n,l,m}=(g, \varphi_{n,l,m})_{L^2(\mathbb{R}^3)}$, and
\begin{equation*}
\mathcal{H}(g)=\sum^{\infty}_{n=0}\sum^{\infty}_{l=0}\sum^{l}_{m=-l}(2n+l+\frac 32)\, g_{n,l,m}\, \varphi_{n, l, m}\, .
\end{equation*}
Using this spectral decomposition, the definition of $(\log(\mathcal{H}+1))^{\frac{2}{s}},\, e^{c(\log(\mathcal{H}+1))^{\frac{2}{s}}},\, e^{c\mathcal{L}}$ is then classical.

\smallskip
\subsection{Triangular effect of the non linear operators}
We study now the algebra property of the nonlinear terms
$$
{\bf \Gamma}(\varphi_{n,0,0},
\varphi_{m,0,0}),
$$
By the same proof of Proposition 2.1 in \cite{GLX_2015}, we have the following triangular effect for the nonlinear Boltzmann operators on the basis $\{\varphi_{n,0,0}\}$.

\begin{proposition}\label{expansion}
The following algebraic identities hold,
\begin{align*}
&(i_1) \quad\,\, {\bf \Gamma}(\varphi_{0,0,0},\varphi_{m,0,0})=
\left(\int^{\frac{\pi}{4}}_0\beta(\theta)
((\cos\theta)^{2m}-1)d\theta\right)\varphi_{m,0,0},\,\, m\in\mathbb{N};\\
&(i_2) \quad\,\, {\bf \Gamma}(\varphi_{n,0,0},\varphi_{0,0,0})=\left(\int^{\frac{\pi}{4}}_0\beta(\theta)((\sin\theta)^{2n}-\delta_{0,n})d\theta\right)\varphi_{n,0,0},\,\,n\in\mathbb{N};\\
&(ii)  \quad\,\, {\bf \Gamma}(\varphi_{n,0,0},\varphi_{m,0,0})=\mu_{n,m}\varphi_{n+m,0,0}, \,\,\text{for}\,\, n\geq1, m\geq1,
\end{align*}
where
\begin{align}\label{mumn}
\mu_{n,m}=\sqrt{\frac{(2n+2m+1)!}{(2n+1)!(2m+1)!}}
\left(\int^{\frac{\pi}{4}}_0\beta(\theta)(\sin\theta)^{2n}(\cos\theta)^{2m}d\theta\right).
\end{align}
\end{proposition}
\begin{remark}\label{lambdasum}
{\color{black}
Obviously, we can deduce from $(i_1)$ and $(i_2)$ of Proposition \ref{expansion} that  }
\begin{equation*}
\quad\forall n\in\mathbb{N}, \quad{\bf \Gamma}(\varphi_{0,0,0},\varphi_{n,0,0})+{\bf \Gamma}(\varphi_{n,0,0},\varphi_{0,0,0})=-\lambda_{n,0}\varphi_{n,0,0}.
\end{equation*}
Where $\lambda_{0,0}=\lambda_{1,0}=0$ and for $n\geq2$,
$$\lambda_{n,0}=\int^{\frac{\pi}{4}}_0\beta(\theta)
(1-(\cos\theta)^{2n}-(\sin\theta)^{2n})d\theta.$$
From Proposition 2.1 in \cite{GL-2015}, there exists a $c_0>0$ dependent only on $s$, such that, for $n\geq2$,
\begin{equation}\label{estilambda}
c_0(\log(2n+\frac{5}{2}))^{\frac{2}{s}}\leq
{\color{black} \lambda_{n,0}}
\leq\frac{1}{c_0}(\log(2n+\frac{5}{2}))^{\frac{2}{s}}.
\end{equation}
This shows that the linearized radially symmetric spatially homogeneous Boltzmann operator with Debye-Yukawa potential was shown to
 behave as a fractional logarithmic harmonic oscilator $\left(\log(\mathcal{H}+1)\right)^{\frac{2}{s}}$.

\end{remark}


\subsection{Explicit solution of the Cauchy problem}\label{subsectionsolution}

Now we solve explicitly the Cauchy problem associated to the non-cutoff radial symmetric spatially
homogeneous Boltzmann equation with Maxwellian molecules for a small $L^2$-initial radial data.

We search a radial solution to the Cauchy problem \eqref{eq-1} in the form
$$
g(t)=\sum^{+\infty}_{n=0}g_{n}(t)\varphi_{n,0,0},
$$
with initial data
$$
g(0)=\sum^{+\infty}_{n=0}\left(g_0,\varphi_{n,0,0}\right)_{L^2(\mathbb{R}^3)}\varphi_{n,0,0}\in L^2(\mathbb{R}^3),
$$
where
$$
g_n(t)=\left(g(t),\, \varphi_{n,0,0}\right)_{L^2(\mathbb{R}^3)}.
$$
It follows from Proposition \ref{expansion} and Remark \ref{lambdasum} that, for convenable radial symmetric function $g$, we have
\begin{align*}
{\bf \Gamma}(g,g)
&=-\sum^{+\infty}_{n=0}
g_0(t)g_n(t)
\lambda_{n,0}\varphi_{n,
0,0}\\
&\quad+\sum^{+\infty}_{n=1}\sum^{+\infty}_{m=1}g_n(t)g_m(t)\mu_{n,m}
\varphi_{m+n,0,0},
\end{align*}
where $\mu_{n,m}$ was defined in \eqref{mumn}.  This implies that,
\begin{align*}
{\bf \Gamma}(g,g)
&=\sum^{+\infty}_{n=0}
\Big[
-g_0(t)g_n(t)\lambda_{n,0}
+\sum_{\substack{k+l=n\\k\geq1,l\geq1}}g_k(t)g_l(t)\mu_{k,l}
\Big]
\varphi_{n,0,0}.
\end{align*}
For radial symmetric function $g$, we also have
\begin{equation*}
\mathcal{L}g=\sum^{+\infty}_{n=0}\lambda_{n,0}\,g_n(t)\,\varphi_{n,0,0}.
\end{equation*}
Formally, we take inner product with $\varphi_{n,0,0}$ on both sides of \eqref{eq-1},
we find that the functions $\{g_{n}(t)\}$ satisfy the following infinite system
of the differential equations
\begin{equation}\label{ODE-1}
\partial_t g_n(t)+\lambda_{n,0}\,g_n(t)=-g_0(t)g_n(t)\lambda_{n,0}
+\sum_{\substack{k+l=n\\k\geq1,l\geq1}}g_k(t)g_l(t)\mu_{k,l},\quad\forall n\in\mathbb{N},
\end{equation}
with initial data
$$g_n(0)=\left(g_0,\varphi_{n,0,0}\right)_{L^2(\mathbb{R}^3)}.$$
Consider that $g_0\in\mathcal{N}^{\perp}$, we have
$$g_0(0)=g_1(0)=0.$$
The infinite system of the differential equations \eqref{ODE-1} reduces to be
\begin{equation}\label{ODE-2}
\left\{ \begin{aligned}
         &g_0(t)=g_1(t)=0,\\
         &\partial_t g_n(t)+\lambda_{n,0}\,g_n(t)=\sum_{\substack{k+l=n\\k\geq1,l\geq1}}g_k(t)g_l(t)\mu_{k,l},\quad\forall n\geq2, \\
         &g_n(0)=\left(g_0,\varphi_{n,0,0}\right)_{L^2(\mathbb{R}^3)}.
\end{aligned} \right.
\end{equation}
On the right hand side of the second equation in \eqref{ODE-2}, the indices $k$ and $l$ are always less than $n$, then this system of the differential equations is triangular, which can be explicitly solved while solving a sequence of linear differential equations.

The proof of Theorem \ref{trick} is reduced to prove the convergence of following series
\begin{equation}\label{ODE}
g(t)=\sum^{+\infty}_{n=2}g_n(t)\varphi_{n, 0,0}
\end{equation}
in the convenable function space.


\section{The sharp trilinear estimates for the radially symmetric Boltzmann operator}\label{S3}
To prove the convergence of the formal solution obtained in the precedent section, we need to estimate the following trilinear terms
$$
\left({\bf \Gamma}(f,g),h\right)_{L^2(\mathbb{R}^3)},
\,\,\,f,g,h\in\mathscr{S}_r(\mathbb{R}^3)\cap\mathcal{N}^{\perp}.
$$
By a proof similar to that in Lemma 3.1 in \cite{NYKC3}, we present the properties of the eigenvalues of the linearized radially symmetric Boltzmann operator $\mathcal{L}$, which is a basic tool in the proof of the trilinear estimate with exponential weighted.
\begin{lemma}\label{eigen}
The eigenvalues of the linearized radially symmetric Boltzmann operator $\mathcal{L}$
$$\lambda_{n,0}=\int_{|\theta|\le \frac{\pi}{4}}\beta(\theta)\Big(1
-(\sin\theta)^{2n}-(\cos\theta)^{2n}\Big)d\theta$$
satisfy to the following estimate
$$\forall k,l\geq2,\quad \lambda_{k,0}+\lambda_{l,0}>\lambda_{k+l,0}.$$
\end{lemma}
\begin{proof}
Since $\beta(\theta)>0$, we only need to prove that, for $\theta\in[-\frac{\pi}{4},\frac{\pi}{4}]\setminus{0}$, $\forall k,l\geq2,$
\begin{equation}\label{inequ}
1+(\cos\theta)^{2k+2l}+(\sin\theta)^{2k+2l}-(\cos\theta)^{2k}-(\sin\theta)^{2k}-(\cos\theta)^{2l}-(\sin\theta)^{2l}>0.
\end{equation}
By a proof similar to that in Lemma 3.1 in \cite{NYKC3}, the estimate \eqref{inequ} follows. This ends the proof of Lemma \ref{eigen}.
\end{proof}

The following lemma is instrumental in the proof of the trilinear estimates:
\begin{lemma}\label{sum}
For $n\geq2$ and $\mu_{k,l}$ was defined in \eqref{mumn} with $k, l\in\mathbb{N}$, $0<s\leq2$, we have
\begin{equation}\label{estimatesum}
\sum_{\substack{k+l=n\\k\geq1,l\geq1}}\frac{|\mu_{k,l}|^2}{(\log(2l+\frac{5}{2}))^{2/s}}\lesssim (\log(2n+\frac{5}{2}))^{2/s}.
\end{equation}
\end{lemma}
We prove this Lemma in Section \ref{section estimate}.

The sharp trilinear estimates for the radially symmetric Boltzmann operator can be derived from the result of Lemma \ref{sum}.
\begin{proposition}\label{estimatetri}
For $0<s\leq2$,
there exists a positive $C>0$, such that for all $f,g,h\in\mathscr{S}_r(\mathbb{R}^3)\cap\mathcal{N}^{\perp}$,
\begin{align*}
|(\Gamma(f,g),h)_{L^2}|\leq C\|f\|_{L^2}\|(\log(\mathcal{H}+1))^{\frac{1}{s}}g\|_{L^2}\|(\log(\mathcal{H}+1))^{\frac{1}{s}}h\|_{L^2},
\end{align*}
and for any $t\geq0$, $n\geq2$,
\begin{align*}
&|(\Gamma(f,g),e^{t\mathcal{L}}\mathbb{S}_nh)_{L^2}|\\
\leq& C\|e^{\frac{t}{2}\mathcal{L}}\mathbb{S}_{n-2}f\|_{L^2}
\|e^{\frac{t}{2}\mathcal{L}}(\log(\mathcal{H}+1))^{\frac{1}{s}}g\mathbb{S}_{n-2}g\|_{L^2}
\|e^{\frac{t}{2}\mathcal{L}}(\log(\mathcal{H}+1))^{\frac{1}{s}}\mathbb{S}_nh\|_{L^2},
\end{align*}
where $\mathcal{L}$ is the linearized non-cutoff Boltzmann operator, $\mathcal{H}=-\triangle +\frac{|v|^2}{4}$ is the 3-dimensional harmonic oscillator
and $\mathbb{S}_n$ is the orthogonal projector onto the $n+1$ energy levels
$$\mathbb{S}_{n}f=\sum^{n}_{k=0}(f,\varphi_{k,0,0})_{L^2}\varphi_{k,0,0},\quad
e^{t\mathcal{L}}\mathbb{S}_nf=\sum^{n}_{k=0}e^{\lambda_{k,0}t}(f,\varphi_{k,0,0})_{L^2}\varphi_{k,0,0}.$$
\end{proposition}
\begin{proof}
Let $f,g,h\in\mathscr{S}_r(\mathbb{R}^3)\cap\mathcal{N}^{\perp}$ be the radial Schwartz functions, by using the spectral decomposition, we obtain
$$f=\sum^{+\infty}_{n=2}(f,\varphi_{n,0,0})_{L^2}\varphi_{n,0,0},\quad g=\sum^{+\infty}_{n=2}(g,\varphi_{n,0,0})_{L^2}\varphi_{n,0,0},\quad h=\sum^{+\infty}_{n=2}(h,\varphi_{n,0,0})_{L^2}\varphi_{n,0,0}.$$
In convenience, we rewrite  $f_n=(f,\varphi_{n,0,0})_{L^2}, g_n=(g,\varphi_{n,0,0})_{L^2},h_n=(h,\varphi_{n,0,0})_{L^2}.$   We can deduce from Proposition \ref{expansion} that,
\begin{align*}
\Gamma(f,g)
&=\sum^{+\infty}_{k=2}\sum^{+\infty}_{l=2}f_kg_l\Gamma(\varphi_{k,0,0},\varphi_{l,0,0})\\
&=\sum^{+\infty}_{k=2}\sum^{+\infty}_{l=2}f_kg_l\mu_{k,l}\varphi_{k+l,0,0}\\
&=\sum^{+\infty}_{n=4}\Big(\sum_{\substack{k+l=n\\k\geq2,l\geq2}}\mu_{k,l}f_kg_l\Big)\varphi_{n,0,0}.
\end{align*}
Applying the orthogonal property of $\varphi_{n,0,0}$, it follows that,
\begin{align*}
(\Gamma(f,g),h)_{L^2}=\sum^{+\infty}_{n=4}\Big(\sum_{\substack{k+l=n\\k\geq2,l\geq2}}\mu_{k,l}f_kg_l\Big)h_n.
\end{align*}
We use the Cauchy-Schwarz inequality,
\begin{align*}
&|(\Gamma(f,g),h)_{L^2}|\\
\leq&\sum^{+\infty}_{k=2}\sum^{+\infty}_{l=2}|\mu_{k,l}||f_k||g_l||h_{k+l}|\\
\leq&\left(\sum^{+\infty}_{l=2}\Big(\log(2l+\frac{5}{2})\Big)^{\frac{2}{s}}|g_l|^2\right)^{\frac{1}{2}}
\left(\sum^{+\infty}_{l=2}\frac{1}{\Big(\log(2l+\frac{5}{2})\Big)^{\frac{2}{s}}}\Big(\sum^{+\infty}_{k=2}|\mu_{k,l}||f_k||h_{k+l}|\Big)^2\right)^{\frac{1}{2}}\\
\leq&\left(\sum^{+\infty}_{l=2}\Big(\log(2l+\frac{5}{2})\Big)^{\frac{2}{s}}|g_l|^2\right)^{\frac{1}{2}}
\left(\sum^{+\infty}_{k=2}|f_k|^2\right)^{\frac{1}{2}}\\
&\qquad\times\left(\sum^{+\infty}_{l=2}\sum^{+\infty}_{k=2}\frac{|\mu_{k,l}|^2}{\Big(\log(2l+\frac{5}{2})\Big)^{\frac{2}{s}}}|h_{k+l}|^2\right)^{\frac{1}{2}}\\
=&\|(\log(\mathcal{H}+1))^{\frac{1}{s}}g\|_{L^2}\|f\|_{L^2}
\left(\sum^{+\infty}_{n=4}\Big(\sum_{\substack{k+l=n\\k\geq2,l\geq2}}\frac{|\mu_{k,l}|^2}{\Big(\log(2l+\frac{5}{2})\Big)^{\frac{2}{s}}}\Big)|h_n|^2\right)^{\frac{1}{2}}.
\end{align*}
It follows from  Lemma \ref{sum} that
$$\sum_{\substack{k+l=n\\k\geq2,l\geq2}}\frac{|\mu_{k,l}|^2}{\Big(\log(2l+\frac{5}{2})\Big)^{\frac{2}{s}}}\lesssim \Big(\log(2n+\frac{5}{2})\Big)^{\frac{2}{s}},$$
then
$$\left(\sum^{+\infty}_{n=4}\Big(\sum_{\substack{k+l=n\\k\geq2,l\geq2}}\frac{|\mu_{k,l}|^2}{\Big(\log(2l+\frac{5}{2})\Big)^{\frac{2}{s}}}\Big)|h_n|^2\right)^{\frac{1}{2}}\lesssim\|(\log(\mathcal{H}+1))^{\frac{1}{s}}h\|_{L^2}.$$
This implies that, there exists $C>0$,
\begin{align*}
|(\Gamma(f,g),h)_{L^2}|\leq C\|(\log(\mathcal{H}+1))^{\frac{1}{s}}g\|_{L^2}\|f\|_{L^2}\|(\log(\mathcal{H}+1))^{\frac{1}{s}}h\|_{L^2}.
\end{align*}
On the other hand, we consider the inequality with exponential weighted and apply the orthogonal property of $\varphi_{n,0,0}$ again that
\begin{align*}
(\Gamma(f,g),e^{t\mathcal{L}}\mathbb{S}_nh)=\sum^{n}_{m=4}e^{\lambda_{m,0}t}\Big(\sum_{\substack{k+l=m\\k\geq2,l\geq2}}\mu_{k,l}f_kg_l\Big)h_m.
\end{align*}
Then
\begin{align*}
&|(\Gamma(f,g),e^{t\mathcal{L}}\mathbb{S}_nh)|\\
\leq&\sum^{n-2}_{l=2}\sum^{n-l}_{k=2}|g_l||\mu_{k,l}||f_k||h_{k+l}|e^{\lambda_{k+l,0}t}\\
\leq&\left(\sum^{n-2}_{l=2}e^{\lambda_{l,0}t}\Big(\log(2l+\frac{5}{2})\Big)^{\frac{2}{s}}|g_l|^2\right)^{\frac{1}{2}}\\
&\times\left(\sum^{n-2}_{l=2}
\frac{ e^{-\lambda_{l,0}t}}{\Big(\log(2l+\frac{5}{2})\Big)^{\frac{2}{s}}}
\Big(\sum^{n-l}_{k=2}e^{\lambda_{k+l,0}t}|\mu_{k,l}||f_k||h_{k+l}|\Big)^2\right)^{\frac{1}{2}}\\
\leq&\left(\sum^{n-2}_{l=2}e^{\lambda_{l,0}t}\Big(\log(2l+\frac{5}{2})\Big)^{\frac{2}{s}}|g_l|^2\right)^{\frac{1}{2}}\left(\sum^{n-2}_{k=2}e^{\lambda_{k,0}t}|f_k|^2\right)^{\frac{1}{2}}\\
&\times\left(\sum^{n-2}_{l=2}
\frac{1}{\Big(\log(2l+\frac{5}{2})\Big)^{\frac{2}{s}}}
\sum^{n-l}_{k=2}e^{2\lambda_{k+l,0}t-\lambda_{l,0}t-\lambda_{k,0}t}|\mu_{k,l}|^2|h_{k+l}|^2\right)^{\frac{1}{2}}.
\end{align*}
Since by Lemma \ref{eigen}, for all $k\geq2,l\geq2,$
$$ \lambda_{k+l,0}-\lambda_{l,0}-\lambda_{k,0}\leq0,$$
one can verify that, for $t\geq0$,
\begin{align*}
&|(\Gamma(f,g),e^{t\mathcal{L}}\mathbb{S}_nh)|\\
\leq &\|e^{\frac{t}{2}\mathcal{L}}\mathbb{S}_{n-2}f\|_{L^2}
\|e^{\frac{t}{2}\mathcal{L}}(\log(\mathcal{H}+1))^{\frac{1}{s}}\mathbb{S}_{n-2}g\|_{L^2}\\
&\qquad\times\left(\sum^{n-2}_{l=2}
\frac{1}{\Big(\log(2l+\frac{5}{2})\Big)^{\frac{2}{s}}}
\sum^{n-l}_{k=2}e^{\lambda_{k+l,0}t}|\mu_{k,l}|^2|h_{k+l}|^2\right)^{\frac{1}{2}}\\
=&\|e^{\frac{t}{2}\mathcal{L}}\mathbb{S}_{n-2}f\|_{L^2}
\|e^{\frac{t}{2}\mathcal{L}}(\log(\mathcal{H}+1))^{\frac{1}{s}}\mathbb{S}_{n-2}g\|_{L^2}\\
&\qquad\times
\left(\sum^{n}_{m=4}e^{\lambda_{m,0}t}\Big(\sum_{\substack{k+l=m\\k\geq2,l\geq2}}\frac{|\mu_{k,l}|^2}{\Big(\log(2l+\frac{5}{2})\Big)^{\frac{2}{s}}}\Big)|h_m|^2\right)^{\frac{1}{2}}.
\end{align*}
By using Lemma \ref{sum} again that, for $m\geq4$,
$$\sum_{\substack{k+l=m\\k\geq2,l\geq2}}\frac{|\mu_{k,l}|^2}{\Big(\log(2l+\frac{5}{2})\Big)^{\frac{2}{s}}}\lesssim \Big(\log(2m+\frac{5}{2})\Big)^{\frac{2}{s}}.$$
We conclude, for $t\geq0$, $n\geq2$,
\begin{align*}
&|(\Gamma(f,g),e^{t\mathcal{L}}\mathbb{S}_nh)_{L^2}|\\
\leq& C\|e^{\frac{t}{2}\mathcal{L}}\mathbb{S}_{n-2}f\|_{L^2}
\|e^{\frac{t}{2}\mathcal{L}}(\log(\mathcal{H}+1))^{\frac{1}{s}}\mathbb{S}_{n-2}g\|_{L^2}
\|e^{\frac{t}{2}\mathcal{L}}(\log(\mathcal{H}+1))^{\frac{1}{s}}\mathbb{S}_nh\|_{L^2}.
\end{align*}
This ends the proof of Proposition \ref{estimatetri}.
\end{proof}

\begin{remark}\label{3.4}
From the remark \ref{estimatesum}, the linearized radially symmetric spatially homogeneous Boltzmann operator with Debye-Yukawa potential $\mathcal{L}$ was shown to
 behave as a fractional logarithmic harmonic oscilator $\left(\log(\mathcal{H}+1)\right)^{\frac{2}{s}}$, one can verify from Proposition \ref{estimatetri} that {\color{black} there exists a constant $C_1$ such that }
 \begin{align*}
|(\Gamma(f,g),e^{t\mathcal{L}}\mathbb{S}_nh)_{L^2}|\leq {\color{black}C_1} \|e^{\frac{t}{2}\mathcal{L}}\mathbb{S}_{n-2}f\|_{L^2}
\|e^{\frac{t}{2}\mathcal{L}}\mathcal{L}^{\frac{1}{2}}\mathbb{S}_{n-2}g\|_{L^2}
\|e^{\frac{t}{2}\mathcal{L}}\mathcal{L}^{\frac{1}{2}}\mathbb{S}_nh\|_{L^2}.
\end{align*}
\end{remark}


\section{The proof of the main Theorem}\label{S4}

In this section, we study the convergence of the formal solutions obtained on Section \ref{S2} with small $L^2$ initial data which end the proof of Theorem  \ref{trick}.

\subsection{The uniform estimate}
Let $\{g_n(t)\}$ be the solution of \eqref{ODE}, for any $2\leq N\in \mathbb{N}$, set
\begin{equation}\label{SN}
\mathbb{S}_Ng(t)=
\sum^{N}_{n=2}g_n(t)\varphi_{n,0,0},
\end{equation}
then $\mathbb{S}_ng(t),  e^{t\mathcal{L}}\mathbb{S}_ng(t) \in\mathscr{S}_r(\mathbb{R}^3)\bigcap\mathcal{N}^{\perp}$,

Multiplying $e^{\lambda_{n,0}t}\overline{g_n}(t)$ on both sides of \eqref{ODE-1} and take summation for $ 2\leq n\leq N$, then Proposition \ref{expansion} and the orthogonality of the basis $\{\varphi_{n,0,0}\}_{n\in \mathbb{N}}$ imply that
\begin{align*}
&\Big(\partial_t (\mathbb{S}_Ng)(t),e^{t\mathcal{L}}\mathbb{S}_N g(t)\Big)_{L^2(\mathbb{R}^3)}+\Big(\mathcal{L}(\mathbb{S}_Ng)(t),e^{t\mathcal{L}}\mathbb{S}_N g(t)\Big)_{L^2(\mathbb{R}^3)}\\
&\qquad=\Big(\Gamma((\mathbb{S}_Ng),(\mathbb{S}_Ng)), e^{t\mathcal{L}}\mathbb{S}_N g(t)\Big)_{L^2(\mathbb{R}^3)}.
\end{align*}
Since $\mathbb{S}_Ng(t)\in\mathscr{S}(\mathbb{R}^3)\bigcap\mathcal{N}^{\perp}$, we have
$$
\Big(\mathcal{L}(\mathbb{S}_Ng)(t), e^{t\mathcal{L}}\mathbb{S}_N g(t)\Big)_{L^2(\mathbb{R}^3)}=
\|e^{\frac{t}{2}\mathcal{L}}\mathcal{L}^{\frac 12}\mathbb{S}_Ng(t)\|^2_{L^2(\mathbb{R}^3)},
$$
we then obtain that
\begin{align*}
\frac{1}{2}\frac{d}{dt}\|e^{\frac{t}{2}\mathcal{L}}\mathbb{S}_Ng(t)\|^2_{L^2}
&+\frac{1}{2}\|e^{\frac{t}{2}\mathcal{L}}\mathcal{L}^{\frac 12}\mathbb{S}_Ng(t)\|^2_{L^2(\mathbb{R}^3)}\\
&=\Big(\Gamma((\mathbb{S}_Ng),(\mathbb{S}_Ng)), e^{t\mathcal{L}}\mathbb{S}_Ng(t)\Big)_{L^2}.
\end{align*}
It follows from Remark \ref{3.4} that, for any $N\geq2$, $t\geq0$,
\begin{align}\label{ele}
\frac{1}{2}\frac{d}{dt}\|e^{\frac{t}{2}\mathcal{L}}\mathbb{S}_Ng\|^2_{L^2}
&+\frac{1}{2}\|e^{\frac{t}{2}\mathcal{L}}\mathcal{L}^{\frac{1}{2}}\mathbb{S}_{N}g\|^2_{L^2}\nonumber\\
&\leq {\color{black} C_1} \|e^{\frac{t}{2}\mathcal{L}}\mathbb{S}_{N-2}g\|_{L^2}
\|e^{\frac{t}{2}\mathcal{L}}\mathcal{L}^{\frac{1}{2}}\mathbb{S}_{N}g\|^2_{L^2}.
\end{align}

\begin{proposition}\label{induction}
There exists $\epsilon_0>0$ such that for all $0<\epsilon\leq\epsilon_0$, ,\, $g_0\in L^2\bigcap\mathcal{N}^{\perp}$ with $\|g_0\|_{L^2}\leq\epsilon$,
\begin{equation*}
\|e^{\frac{t}{2}\mathcal{L}}\mathbb{S}_{N}g(t)\|^2_{L^2(\mathbb{R}^3)}
+\frac{1}{2}\int^t_0\|e^{\frac{t}{2}\mathcal{L}}\mathcal{L}^{\frac{1}{2}}
\mathbb{S}_{N}g(\tau)\|^2_{L^2}d\tau
\le\|g_0\|^2_{L^2(\mathbb{R}^3)},
\end{equation*}
for any $t\geq0,\, N\geq 0$.
\end{proposition}

\begin{proof} We prove the Proposition by induction on $N$.

{\bf 1). For $N\le 2$.} we have
$\|e^{\frac{t}{2}\mathcal{L}}\mathbb{S}_{0}g\|^2_{L^2}=|g_0(t)|^2=0,\,$
$\|e^{\frac{t}{2}\mathcal{L}}\mathbb{S}_{1}g\|^2_{L^2}=
|g_0(t)|^2+|g_1(t)|^2=0,$
and
$$
\|e^{\frac{t}{2}\mathcal{L}}\mathbb{S}_{2}g\|^2_{L^2}=e^{\lambda_{0,2}t}|g_2(t)|^2
=e^{-\lambda_{0,2}t}|g_2(0)|^2\leq|\left(g_0,\varphi_{2,0,0}\right)_{L^2(\mathbb{R}^3)}|^2\leq\|g_0\|^2_{L^2(\mathbb{R}^3)}.
$$

{\bf 2). For $N> 2$.} We want to prove that
$$
\|e^{\frac{t}{2}\mathcal{L}}\mathbb{S}_{N-1}g\|_{L^2}\leq\epsilon\leq\epsilon_0,
$$
imply
$$
\|e^{\frac{t}{2}\mathcal{L}}\mathbb{S}_{N}g\|_{L^2}\leq\epsilon.
$$
Take now $\epsilon_0>0$ such that,
$$
0<\epsilon_0\leq\frac{1}{4 {\color{black} C_1}}
$$
{\color{black} where $ C_1$ is defined in Remark \ref{3.4}.}
Then we deduce from \eqref{ele} that
\begin{align*}
\frac{1}{2}\frac{d}{dt}\|e^{\frac{t}{2}\mathcal{L}}\mathbb{S}_Ng(t)\|^2_{L^2}
&+\frac{1}{2}\|e^{\frac{t}{2}\mathcal{L}}\mathcal{L}^{\frac{1}{2}}\mathbb{S}_{N}g\|^2_{L^2}\nonumber\\
&\leq {\color{black} C_1}\|e^{\frac{t}{2}\mathcal{L}}\mathbb{S}_{N-2}g\|_{L^2}
\|e^{\frac{t}{2}\mathcal{L}}\mathcal{L}^{\frac{1}{2}}\mathbb{S}_{N}g\|^2_{L^2}\\
&\leq\frac{1}{4}\|e^{\frac{t}{2}\mathcal{L}}
\mathcal{L}^{\frac{1}{2}}\mathbb{S}_{N}g\|^2_{L^2},
\end{align*}
therefore,
\begin{equation}\label{indu1}
\frac{d}{dt}\|e^{\frac{t}{2}\mathcal{L}}\mathbb{S}_Ng(t)\|^2_{L^2}
+\frac{1}{2}\|e^{\frac{t}{2}\mathcal{L}}\mathcal{L}^{\frac{1}{2}}\mathbb{S}_{N}g\|^2_{L^2}\le 0.
\end{equation}
This ends the proof of the Proposition \ref{induction}.
\end{proof}

\subsection{Existence of the weak solution}\,\,
We prove now the convergence of the sequence
$$
g(t)=\sum^{+\infty}_{n=2} g_{n}(t)\varphi_{n,0,0}
$$
defined in \eqref{ODE}.  For all $N\geq2$,
$\mathbb{S}_{N}g(t)$ satisfies the following Cauchy problem
\begin{equation}\label{eq-2}
\left\{ \begin{aligned}
         &\partial_t \mathbb{S}_{N}g+\mathcal{L}(\mathbb{S}_{N}g)=\mathbb{S}_{N}\Gamma(\mathbb{S}_{N}g, \mathbb{S}_{N}g),\,\\
         &\mathbb{S}_{N}g(0)=\sum^{N}_{n=2}\left(g_0,\varphi_{n,0,0}\right)_{L^2(\mathbb{R}^3)}\varphi_{n,0,0}.
\end{aligned} \right.
\end{equation}
By Proposition \ref{induction} and the orthogonality of the basis $(\varphi_{n,0,0})_{n\in\mathbb{N}}$, for all $t>0$,
\begin{equation*}
\sum^{N}_{n=2}e^{\lambda_{n,0}t}|g_n(t)|^2
+\frac{1}{2}\int^t_0\left(\sum^{N}_{n=2}e^{\lambda_{n,0}t}\lambda_{n,0}|g(\tau)|^2\right)d\tau
\le\|g_0\|^2_{L^2(\mathbb{R}^3)}.
\end{equation*}
Since $\lambda_{n,0}>0$, for all $t\geq0$, we have
\begin{equation}\label{Lestimate}
\sum^{N}_{n=2}|g_n(t)|^2
+\frac{1}{2}\int^t_0\left(\sum^{N}_{n=2}\lambda_{n,0}|g(\tau)|^2\right)d\tau
\le\|g_0\|^2_{L^2(\mathbb{R}^3)}.
\end{equation}
The orthogonality of the basis $(\varphi_{n,0,0})_{n\in\mathbb{N}}$ implies that
$$
\|\mathbb{S}_{N}g(t)\|^2_{L^2(\mathbb{R}^3)}=
\sum^{N}_{n=2}|g_n(t)|^2,\quad\|\mathcal{L}^{\frac{1}{2}}\mathbb{S}_{N}g(t)\|^2_{L^2(\mathbb{R}^3)}=
\sum^{N}_{n=2}\lambda_{n,0}|g_n(t)|^2.
$$
By using the monotone convergence theorem, the sequence
$$
g(t)=\sum^{+\infty}_{n=2}\,\,
g_n(t)\varphi_{n,0,0}
$$
is convergent and for any $t\ge 0$,
$$
\lim_{N\to \infty}\|\mathbb{S}_{N}g-g\|_{L^\infty(]0, t[; L^2(\mathbb{R}^3))}=0
$$
and
$$
\lim_{N\to \infty}\|\mathcal{L}^{\frac{1}{2}}(\mathbb{S}_{N}g-g)\|_{L^2(]0, t[; L^2(\mathbb{R}^3))}=0.
$$
%
%
For any $\phi(t)\in\,C^1\Big(\mathbb{R}_+,\mathscr{S}(\mathbb{R}^3)\Big)$,  we have
\begin{align*}
&\Big|\int^t_0\Big(\mathbb{S}_{N}\Gamma(\mathbb{S}_{N}g, \mathbb{S}_{N}g)
-\Gamma(g, g),
\phi(\tau)\Big)_{L^2(\mathbb{R}^3)}d\tau\Big|\\
&\le\Big|\int^t_0\Big(\Gamma(\mathbb{S}_{N}g, \mathbb{S}_{N}g),\mathbb{S}_{N}\phi(\tau)-\phi(\tau)\Big)_{L^2}d\tau\Big|\\
&\qquad+\Big|\int^t_0\Big(\Gamma(\mathbb{S}_{N}g-g, \mathbb{S}_{N}g),\phi(\tau)\Big)_{L^2}d\tau\Big|+\Big|\int^t_0\Big(\Gamma(g, \mathbb{S}_{N}g-g),\phi(\tau)\Big)_{L^2}d\tau\Big|.
\end{align*}
By the estimate \eqref{Lestimate} and the orthogonality of the basis $(\varphi_{n,0,0})_{n\in\mathbb{N}}$, one can verify that,
\begin{align*}
&\Big|\int^t_0\Big(\mathbb{S}_{N}\Gamma(\mathbb{S}_{N}g, \mathbb{S}_{N}g)-\Gamma(g, g),
\phi(\tau)\Big)_{L^2(\mathbb{R}^3)}d\tau\Big|\\
&\leq\,C\int^t_0\|\mathbb{S}_{N}g\|_{L^2}\|\mathcal{L}^{\frac 12}\mathbb{S}_{N}g\|_{L^2} \Big\|\mathcal{L}^{\frac 12}(
\mathbb{S}_{N}\phi-\phi)\Big\|_{L^2(\mathbb{R}^3)} dt\\
&\qquad+C\int^t_0 \|\mathbb{S}_{N}g-g\|_{L^2}\|\mathcal{L}^{\frac 12}\mathbb{S}_{N}g\|_{L^2} \|\mathcal{L}^{\frac 12}\phi\|_{L^2(\mathbb{R}^3)}
dt\\
&\qquad+C\int^t_0\|g\|_{L^2}\|\mathcal{L}^{\frac 12}(\mathbb{S}_{N}g-g)\|_{L^2(\mathbb{R}^3)}\|\mathcal{L}^{\frac 12}\phi\|_{L^2(\mathbb{R}^3)}
dt\\
&\leq\,C\|g_0\|^2_{L^2} \Big\|\mathcal{L}^{\frac 12}(
\mathbb{S}_{N}\phi-\phi)\Big\|_{L^2(]0, t[; L^2(\mathbb{R}^3))} \\
&\qquad+C\|\mathbb{S}_{N}g-g\|_{L^\infty(]0, t[; L^2)}\|g_0\|_{L^2} \|\mathcal{L}^{\frac 12}\phi\|_{L^2(]0, t[; L^2)}\\
&\qquad+C\|g_0\|_{L^2}\|\mathcal{L}^{\frac 12}(\mathbb{S}_{N}g-g)\|_{L^2(]0, t[; L^2)}\|\mathcal{L}^{\frac 12}\phi\|_{L^2(]0, t[; L^2)}.
\end{align*}
Let $N\rightarrow+\infty$ in \eqref{eq-2}, we conclude that, for any $\phi(t)\in\,C^1\Big(\mathbb{R}_+,\mathscr{S}(\mathbb{R}^3)\Big)$,
\begin{align*}
&\Big(g(t), \phi(t)\Big)_{L^2(\mathbb{R}^3)}-\Big(g(0), \phi(0)\Big)_{L^2(\mathbb{R}^3)}\nonumber\\
&=-\int^t_{0}\Big(\mathcal{L}g(\tau), \phi(\tau)\Big)_{L^2(\mathbb{R}^3)}d\tau+\int^t_{0}\Big(\Gamma(g(\tau),g(\tau)), \phi(\tau)\Big)_{L^2(\mathbb{R}^3)}d\tau,
\end{align*}
which shows $g\in L^\infty(]0, +\infty[; L^2(\mathbb{R}^3))$ is a global weak solution of Cauchy problem \eqref{eq-1}.

\subsection{Regularity of the solution.}\,\, For $\mathbb{S}_{N}g$ defined in \eqref{SN},\,since
$$
\lambda_{n,0}\geq\lambda_{2,0}>0, \, \forall\,n\geq2,
$$
we deduce from the formulas  \eqref{indu1} and the orthogonality of the basis $(\varphi_{n,0,0})_{n\in\mathbb{N}}$ that
\begin{align*}
&\frac{d}{dt}\|e^{\frac{t}{2}\mathcal{L}}\mathbb{S}_Ng(t)\|^2_{L^2}
+\frac{\lambda_{2,0}}{2}\|e^{\frac{t}{2}\mathcal{L}}\mathbb{S}_{N}g\|^2_{L^2}\\
&\leq\frac{d}{dt}\|e^{\frac{t}{2}\mathcal{L}}\mathbb{S}_Ng(t)\|^2_{L^2}
+\frac{1}{2}\sum^{N}_{n=2}e^{\lambda_{n,0}t}\lambda_{n,0}|g_n(t)|^2\\
&= \frac{d}{dt}\|e^{\frac{t}{2}\mathcal{L}}\mathbb{S}_Ng(t)\|^2_{L^2}
+\frac{1}{2}\|e^{\frac{t}{2}\mathcal{L}}\mathcal{L}^{\frac{1}{2}}\mathbb{S}_{N}g\|^2_{L^2}\leq0.
\end{align*}
We have then
$$
\frac{d}{dt}\Big(e^{\frac{\lambda_{2,0}t}{2}}\|e^{\frac{t}{2}\mathcal{L}}
\mathbb{S}_Ng(t)\|^2_{L^2}\Big)
\leq 0,
$$
it deduces that for any $t>0$, and $N\in\mathbb{N}$,
\begin{equation*}
\|e^{\frac{t}{2}\mathcal{L}}\mathbb{S}_{N}g(t)\|_{L^2(\mathbb{R}^3)}
\le\,e^{-\frac{\lambda_{2,0}t}{4}}\|g_0\|_{L^2(\mathbb{R}^3)}.
\end{equation*}
By using the monotone convergence theorem and the formula \eqref{estilambda}, we conclude that, there exists a constant $c_0>0$, such that
\begin{equation}\label{final}
\|e^{c_0t(\log(\mathcal{H}+1))^{\frac{2}{s}}}g(t)\|_{L^2(\mathbb{R}^3)}\le\,
e^{-\frac{\lambda_{2,0}t}{4}}\|g_0\|_{L^2(\mathbb{R}^3)}.
\end{equation}
This is the formula \eqref{rate}.

For the case $1)$ when $0<s\leq2$, the orthogonality of the basis $(\varphi_{n,0,0})_{n\in\mathbb{N}}$ implies that,
\begin{align*}
\|g\|_{Q^{2c_0t}(\mathbb{R}^3)}&=\|\left(\mathcal{H}+1\right)^{c_0t}g\|_{L^2(\mathbb{R}^3)}
=\sum^{+\infty}_{n=0}(2n+\frac{5}{2})^{c_0t}|g_n|^2\\
&\leq \sum^{+\infty}_{n=0}e^{c_0t(\log(2n+\frac{5}{2}))^{\frac{2}{s}}}|g_n|^2
=\|e^{c_0t(\log(\mathcal{H}+1))^{\frac{2}{s}}}g(t)\|_{L^2(\mathbb{R}^3)}\\
&\le\,
e^{-\frac{\lambda_{2,0}t}{4}}\|g_0\|_{L^2(\mathbb{R}^3)}.
\end{align*}
This is the formula \eqref{rate1}.

For the part $2)$ of Theorem \ref{trick}, in case $0<s<2$,
we deduce from Proposition \ref{sobolev-type} and formula \eqref{final},
the formula \eqref{rate2} follows.

The proof of Theorem \ref{trick} is completed.




\section{Estimate on the nonnear eigenvalue}\label{section estimate}

In this section, we provide the proof of Lemma \ref{sum}.
\begin{lemma}\label{sum-b}
For $n\geq2$ and $\mu_{k,l}$ was defined in \eqref{mumn}
with $k, l\in\mathbb{N}$, $0<s\leq2$, we have
\begin{equation}\label{estimatesum-b}
\sum_{\substack{k+l=n\\k\geq1,l\geq1}}\frac{|\mu_{k,l}|^2}{(\log(2l+\frac{5}{2}))^{2/s}}\lesssim (\log(2n+\frac{5}{2}))^{2/s}.
\end{equation}
\end{lemma}

\begin{proof}

Recall from \eqref{mumn} that
$$
\mu_{k,l}=\sqrt{\frac{(2k+2l+1)!}{(2k+1)!(2l+1)!}}
\left(\int^{\frac{\pi}{4}}_0\beta(\theta)(\sin\theta)^{2k}(\cos\theta)^{2l}d\theta\right),$$
and Beta function \eqref{beta} that $\beta(\theta)\sim(\sin\theta)^{-1}(\log(\sin\theta)^{-1})^{\frac{2}{s}-1},$ we obtain
\begin{equation*}
|\mu_{k,l}|^2\sim\frac{(2k+2l+1)!}{(2k+1)!(2l+1)!}\left(\int^{\frac{\pi}{4}}_0(\log(\sin\theta)^{-1})^{\frac{2}{s}-1}(\sin\theta)^{2k-1}(\cos\theta)^{2l}d\theta\right)^2.
\end{equation*}
By using the substitution rule with $x=(\sin\theta)^2$, then
\begin{align*}
&\int^{\frac{\pi}{4}}_0(\log(\sin\theta)^{-1})^{\frac{2}{s}-1}(\sin\theta)^{2k-1}(\cos\theta)^{2l}d\theta\\
=&2^{-\frac{2}{s}}\int^{\frac{1}{2}}_0\left(\log x^{-1}\right)^{\frac{2}{s}-1}x^{k-1}(1-x)^{l-\frac{1}{2}}dx.
\end{align*}
This shows that,
\begin{equation*}
|\mu_{k,l}|^2\sim\frac{(2k+2l+1)!}{(2k+1)!(2l+1)!}\left(\int^{\frac{1}{2}}_0\left(\log\frac{1}{x}\right)^{\frac{2}{s}-1}x^{k-1}(1-x)^{l-\frac{1}{2}}dx\right)^2.
\end{equation*}
Without loss of generality, we assume $n\gg1$.
We can divide the summation into two parts:
\begin{align}\label{divide}
&\sum_{\substack{k+l=n\\k\geq1,l\geq1}}\frac{|\mu_{k,l}|^2}{(\log(2l+\frac{5}{2}))^{2/s}}
=\sum^{n-1}_{l=1}\frac{|\mu_{n-l,l}|^2}{(\log(2l+\frac{5}{2}))^{2/s}}\nonumber\\
&=\frac{|\mu_{n-1,1}|^2}{(\log\frac{9}{2})^{2/s}}+\sum^{n-1}_{l=2}\frac{|\mu_{n-l,l}|^2}{(\log(2l+\frac{5}{2}))^{2/s}}=\mathbf{H}+\mathbf{I}.
\end{align}
For the estimate of $\mathbf{H}$, since $1-x\sim1$, we have
\begin{align*}
|\mu_{n-1,1}|^2\sim&\frac{(2n+1)!}{(2n-1)!3!}\left(\int^{\frac{1}{2}}_0\left(\log\frac{1}{x}\right)^{\frac{2}{s}-1}x^{n-2}dx\right)^2\\
=&\frac{(2n+1)!}{(2n-1)!3!}\left(\int^{+\infty}_{\log2}x^{\frac{2}{s}-1}e^{-(n-1)x}dx\right)^2\\
=&\frac{(2n+1)!}{(2n-1)!3!}\frac{1}{(n-1)^{4/s}}\left(\int^{+\infty}_{(n-1)\log2}
  u^{\frac{2}{s}-1}e^{-u}du\right)^2.
\end{align*}
It follows from the condition $0<s\leq2$ that,
\begin{equation}\label{11}
\mathbf{H}=\frac{|\mu_{n-1,1}|^2}{(\log\frac{9}{2})^{2/s}}\lesssim\frac{2(2n+1)n}{n^{4/s}4^{n-1}}\lesssim 1.
\end{equation}
We now estimate $\mathbf{I}$.
By using the inequality
$$\left(\int^{\frac{1}{l}}_0+\int^{\frac{1}{2}}_{\frac{1}{l}}\right)^2\leq 2\left(\int^{\frac{1}{l}}_0\right)^2+2\left(\int^{\frac{1}{2}}_{\frac{1}{l}}\right)^2,$$
one can verify that
\begin{align}\label{12}
\mathbf{I}
=&\sum^{n-1}_{l=2}\frac{1}{(\log(2l+\frac{5}{2}))^{2/s}}\frac{(2n+1)!}{(2n-2l+1)!(2l+1)!}\left(\int^{\frac{1}{2}}_0\left(\log x^{-1}\right)^{\frac{2}{s}-1}x^{n-l-1}(1-x)^{l}dx\right)^2\nonumber\\
\leq&\sum^{n-1}_{l=2}\frac{2}{(\log(2l+\frac{5}{2}))^{2/s}}\frac{(2n+1)!}{(2n-2l+1)!(2l+1)!}\Big(\int^{\frac{1}{2}}_{\frac{1}{l}}\left(\log x^{-1}\right)^{\frac{2}{s}-1}x^{n-l-1}(1-x)^ldx\Big)^2\nonumber\\
&+\sum^{n-1}_{l=2}\frac{2}{(\log(2l+\frac{5}{2}))^{2/s}}\frac{(2n+1)!}{(2n-2l+1)!(2l+1)!}\Big(\int^{\frac{1}{l}}_0\left(\log x^{-1}\right)^{\frac{2}{s}-1}x^{n-l-1}(1-x)^ldx\Big)^2\nonumber\\
=&\mathbf{I}_{1}+\mathbf{I}_{2}.
\end{align}
For the part $\mathbf{I}_{1}$.  It is obviously that
\begin{align*}
&\left(\int^{\frac{1}{2}}_{\frac{1}{l}}\left(\log x^{-1}\right)^{\frac{2}{s}-1}x^{n-l-1}(1-x)^ldx\right)^2\\
\lesssim&\left(\log l\right)^{4/s}\left(\int^{\frac{1}{2}}_{\frac{1}{l}}\left(\log(1/x)\right)^{-1}x^{n-l-1}(1-x)^{l+\frac{1}{2}}dx\right)^2\\
\leq&\left(\log l\right)^{4/s}\left(\int^{\frac{1}{2}}_{\frac{1}{l}}\left(\log(1/x)\right)^{-2}x^{-1}dx\right)\left(\int^1_0x^{2n-2l-1}(1-x)^{2l+1}dx\right)\\
=&\frac{(2n-2l-1)!(2l+1)!(\log l)^{4/s}}{(2n+1)!}\left(\int^{\log l}_{\log 2}x^{-2}dx\right)\\
\lesssim&\frac{(2n-2l-1)!(2l+1)!(\log l)^{4/s}}{(2n+1)!}.
\end{align*}
We obtain
\begin{align*}
\mathbf{I}_{1}\lesssim&\sum^{n-1}_{l=2}\frac{1}{(\log(2l+\frac{5}{2}))^{2/s}}\frac{(2n+1)!}{(2n-2l+1)!(2l+1)!}\frac{(2n-2l-1)!(2l+1)!(\log l)^{4/s}}{(2n+1)!}\\
\leq&\left(\log (2n+\frac{5}{2})\right)^{2/s}\sum^{n-1}_{l=2}\frac{1}{(2n-2l+1)(2n-2l)}\\
=&\left(\log (2n+\frac{5}{2})\right)^{2/s}\sum^{n-2}_{k=1}\frac{1}{(2k+1)(2k)}\\
\leq&\left(\log (2n+\frac{5}{2})\right)^{2/s}.
\end{align*}
For the part $\mathbf{I}_{2}$, we divide the summation into two parts,
\begin{align*}
\mathbf{I}_{2}=&\sum^{\leq \frac{n}{2}}_{l=2}\frac{1}{(\log(2l+\frac{5}{2}))^{2/s}}\frac{(2n+1)!}{(2n-2l+1)!(2l+1)!}\Big(\int^{\frac{1}{l}}_0\left(\log x^{-1}\right)^{\frac{2}{s}-1}x^{n-l-1}(1-x)^ldx\Big)^2\\
+&\sum_{\frac{n}{2}<l\leq n-1}\frac{1}{(\log(2l+\frac{5}{2}))^{2/s}}\frac{(2n+1)!}{(2n-2l+1)!(2l+1)!}\Big(\int^{\frac{1}{l}}_0\left(\log x^{-1}\right)^{\frac{2}{s}-1}x^{n-l-1}(1-x)^ldx\Big)^2\\
=&\mathbf{I}_{21}+\mathbf{I}_{22}.
\end{align*}
For the estimation of $\mathbf{I}_{21}$: Transforming $x=\frac{u}{l}$, we have
\begin{align}\label{1}
\int^{\frac{1}{l}}_0\left(\log x^{-1}\right)^{\frac{2}{s}-1}x^{n-l-1}(1-x)^ldx=\int^{1}_0\left(\log \frac{l}{u}\right)^{\frac{2}{s}-1}\frac{u^{n-l-1}}{l^{n-l}}(1-\frac{u}{l})^ldu.
\end{align}
Since
\begin{equation*}
\frac{1}{l^{n-l}}
\left\{
\begin{aligned}
         &=4\frac{1}{2^n},\,\,l=2;\\
         &\leq 9\frac{1}{2^n},\,\,l=3;\\
         &\leq\frac{1}{4^{\frac{n}{2}}}=\frac{1}{2^n},\,\,4\leq l\leq \frac{n}{2},
\end{aligned} \right.
\end{equation*}
we conclude that, for $2\leq l\leq \frac{n}{2}$, $$\frac{1}{l^{n-l}}\lesssim \frac{1}{2^n}.$$
At the same time, since $(1+\frac{1}{l})^{l}$ is increasing with respect to $l$, we obtain, for $0<u<1$,
$$1>(1-\frac{u}{l})^l>(1-\frac{1}{l})^l=\frac{1-\frac{1}{l}}{(1+\frac{1}{l-1})^{l-1}}\geq\frac{1}{4}.$$
It follows from the estimation \eqref{1} that
\begin{align*}
&\int^{\frac{1}{l}}_0\left(\log x^{-1}\right)^{\frac{2}{s}-1}x^{n-l-1}(1-x)^ldx\\
\lesssim &\frac{1}{2^n}\int^{1}_0\left(\log \frac{l}{u}\right)^{\frac{2}{s}-1}u^{n-l-1}du\lesssim\frac{1}{2^n}\left(\log l\right)^{\frac{2}{s}-1}.
\end{align*}
Therefore, we can estimate
\begin{align*}
\mathbf{I}_{21}&\lesssim \sum^{\leq \frac{n}{2}}_{l=2}\frac{1}{(\log(2l+\frac{5}{2}))^{2/s}}\frac{(2n+1)!}{(2n-2l+1)!(2l+1)!}\left(\frac{1}{2^n}\left(\log l\right)^{\frac{2}{s}-1}\right)^2\\
&\lesssim(\log(2n+\frac{5}{2}))^{2/s}\sum^{\leq \frac{n}{2}}_{l=2}\frac{(2n)!}{(2n-2l)!(2l)!}\frac{1}{2^{2n}}\leq (\log(2n+\frac{5}{2}))^{2/s}.
\end{align*}
For the part $\mathbf{I}_{22}$,
\begin{align*}
\mathbf{I}_{22}
\lesssim&\frac{1}{(\log(2n+\frac{5}{2}))^{2/s}}\sum_{\frac{n}{2}<l\leq n-1}\frac{(2n)!}{(2n-2l)!(2l)!}\Big(\int^{\frac{n}{l}}_0\left(\log \frac{n}{u}\right)^{\frac{2}{s}-1}\frac{u^{n-l-1}}{n^{n-l}}(1-\frac{u}{n})^ldu\Big)^2\\
\leq&\frac{1}{(\log(2n+\frac{5}{2}))^{2/s}}\left(\int^2_0\left(\log \frac{n}{u}\right)^{\frac{4}{s}-2}du\right)\\
&\times\left(\int^2_0\sum_{\frac{n}{2}<l\leq n-1}\frac{(2n)!}{(2n-2l)!(2l)!}\frac{u^{2n-2l}}{n^{2n-2l}}(1-\frac{u}{n})^{2l}\frac{du}{u^2}\right).
\end{align*}
For the first integral,
\begin{align*}
\int^2_0\left(\log \frac{n}{u}\right)^{\frac{4}{s}-2}du
&\lesssim\int^2_0(\log n )^{\frac{4}{s}-2}du+\int^2_0\left(\log\frac{1}{u}\right)^{\frac{4}{s}-2}du\\
&\lesssim(\log n )^{\frac{4}{s}-2}+\int^{+\infty}_{\log2}x^{\frac{4}{s}-2}e^{-x}dx\lesssim(\log n )^{\frac{4}{s}-2}.
\end{align*}
We now consider the second integral.  Since for $0<u<2$,
\begin{align*}
\sum_{\frac{n}{2}<l\leq n-1}\frac{(2n)!}{(2n-2l)!(2l)!}\frac{u^{2n-2l}}{n^{2n-2l}}(1-\frac{u}{n})^{2l}&\leq\sum^{2n-2}_{m=0}\frac{(2n)!}{(2n-m)!m!}\frac{u^{2n-m}}{n^{2n-m}}(1-\frac{u}{n})^{m}\\
&=1-(1-\frac{u}{n})^{2n}-2n\frac{u}{n}(1-\frac{u}{n})^{2n-1}.
\end{align*}
Using the Taylor formula for $f(x)=1-(1-x)^{2n}-2nx(1-x)^{2n-1}$ with $0<x<1$,
$$f(x)=f(0)+f'(\theta x)x=2n(2n-1)\theta(1-\theta x)^{2n-2}x^2,  \text{for\,some}\,\,0<\theta<1,$$
we obtain, for $0<u<2$, $n\geq2$,
\begin{align*}
\sum_{\frac{n}{2}<l\leq n-1}\frac{(2n)!}{(2n-2l)!(2l)!}\frac{u^{2n-2l}}{n^{2n-2l}}(1-\frac{u}{n})^{2l}\leq 4u^2.
\end{align*}
This shows that
$$\left(\int^2_0\sum_{\frac{n}{2}<l\leq n-1}\frac{(2n)!}{(2n-2l)!(2l)!}\frac{u^{2n-2l}}{n^{2n-2l}}(1-\frac{u}{n})^{2l}\frac{du}{u^2}\right)\leq 8.$$
Therefore, we have
\begin{align*}
\mathbf{I}_{22}\lesssim\left(\log (2n+\frac{5}{2}) \right)^{\frac{2}{s}}.
\end{align*}
Combining with the estimate of $\mathbf{I}_{1}$, $\mathbf{I}_{21}$, $\mathbf{I}_{22}$ and substituting it into \eqref{12}, we have
$$\mathbf{I}=\sum^{n-1}_{l=2}\frac{|\mu_{n-1,l}|^2}{(\log(2l+\frac{5}{2}))^{2/s}}\leq (\log(2n+\frac{5}{2}))^{2/s}.$$
Reminding the estimate of $\mathbf{H}$ in \eqref{11}, the formula \eqref{estimatesum-b} follows.  This end the proof of Lemma \ref{sum-b}.
\end{proof}


\section{Appendix}\label{appendix}

The important known results but really needed for this paper are presented in this section. For the self-content of paper, we will present some proof of those properties.

\subsection{Shubin spaces}\label{Ap1} We refer the reader to the works \cite{GPR}, \cite{Shubin} for the Shubin spaces.
Let $\tau\in\mathbb{R}$, The Shubin spaces $Q^{\tau}(\mathbb{R}^3)$
can be also characterized through the decomposition
into the Hermite basis:
\begin{align*}
f\in Q^{\tau}(\mathbb{R}^3)
&\Leftrightarrow\,f\in\,L^2(\mathbb{R}^3),
\Bigl\|\Bigl(\mathcal{H}+1\Bigr)^{\frac{\tau}{2}} \, f\Bigr\|_{L^2}<+\infty;
\\
&\Leftrightarrow\, f\in\,L^2(\mathbb{R}^3),\,\,
\Big\|\Big(
(|\alpha|+\frac{5}{2})^{\tau/2} (f,\,H_{\alpha})_{L^2}
\Big)_{\alpha\in\mathbb{N}^3}\Big\|_{l^2}<+\infty,
\end{align*}
where $|\alpha|=\alpha_1+\alpha_2+\alpha_3$,
$$H_{\alpha}(v)=H_{\alpha_1}(v_1)H_{\alpha_2}(v_2)H_{\alpha_3}(v_3),\,\,\alpha\in\mathbb{N}^3,$$
and for $x\in\mathbb{R}$, $n\in \mathbb{N}$,
$$H_{n}(x)=\frac{1}{(2\pi)^{\frac{1}{4}}}\frac{1}{\sqrt{n!}}\Big(\frac{x}{2}-\frac{d}{dx}\Big)^n(e^{-\frac{x^2}{4}}).$$
The following proof is based on the Appendix in \cite{MPX}.
\begin{proof}
Setting $$A_{\pm,j}=\frac{v_j}{2}\mp\frac{d}{dv_j},\quad j=1,2,3,$$
we have, for $\alpha\in\mathbb{N}^3, v\in\mathbb{R}^3,$
$$H_{\alpha}(v)=\frac{1}{\sqrt{\alpha_1!\alpha_2!\alpha_3!}}A^{\alpha_1}_{+,1}H_0(v_1)A^{\alpha_2}_{+,2}H_0(v_2)A^{\alpha_3}_{+,3}H_0(v_3),$$
and for j=1,2,3,
$$A_{+,j}H_{\alpha}=\sqrt{\alpha_j+1}H_{\alpha+e_j},\quad A_{-,j}H_{\alpha}=\sqrt{\alpha_j}H_{\alpha-e_j}(=0\,\text{if}\,\alpha_j=0)$$
where $(e_1,e_2,e_3)$ stands for the canonical basis of $\mathbb{R}^3$.
For the harmonic oscillator $\mathcal{H}=-\triangle +\frac{|v|^2}{4}$ of 3-dimension and $s>0$, one can verify that,
$$\mathcal{H}=\frac{1}{2}\sum^3_{j=1}(A_{+,j}A_{-,j}+A_{-,j}A_{+,j}).$$
Therefore, we have
\begin{align*}
\mathcal{H}H_{\alpha}
&=\frac{1}{2}\sum^3_{j=1}(A_{+,j}A_{-,j}+A_{-,j}A_{+,j})H_{\alpha}\\
&=\frac{1}{2}[\sum^3_{j=1}\sqrt{\alpha_j}A_{+,j}H_{\alpha-e_j}+\sum^3_{j=1}\sqrt{\alpha_j+1}A_{-,j}H_{\alpha+e_j}]\\
&=\frac{1}{2}\sum^3_{j=1}(2\alpha_j+1)H_{\alpha}=\sum^3_{j=1}(\alpha_j+\frac{1}{2})H_{\alpha}.
\end{align*}
By using this spectral decomposition, we conclude that
$$
(\mathcal{H}+1)^{\frac{\tau}{2}} H_{\alpha} = (\lambda_{\alpha}+1)^{\frac{\tau}{2}}H_{\alpha},\,\, \lambda_{\alpha}=\sum^3_{j=1}(\alpha_j+\frac{1}{2}),\,\alpha\in\mathbb{N}^3.
$$
This ends the proof of the another definition of the Shubin space.
\end{proof}

\subsection{Smoothing effects}
Concerning the Shubin spaces introduced in part \ref{Ap1}, we have the following property :
\begin{proposition}\label{sobolev-type}
Let $0<s<2$ and $\tau>0$.
There exists a constant $C=C_s$ such that,
\[
\forall k\geq 1,\quad
\|f\|_{Q^{k}(\mathbb{R}^3)}=\Bigl\|\Bigl(\mathcal{H}+1\Bigr)^{\frac{k}{2}} \, f\Bigr\|_{L^2}\
\leq
  \, e^{C \, \big(\frac{1}{\tau}\big)^{\frac{s}{2-s}}  \, k^{\frac{2}{2-s}}}
  \, \|e^{\tau(\log(\mathcal{H}+1))^{\frac{2}{s}}}f\|_{L^2(\mathbb{R}^3)}
\]
where $\mathcal{H}=-\Delta+{\textstyle \frac{|v|^2}{4}}.$
\end{proposition}

\begin{proof}
Expanding $f$ in the Hermite basis, and noting $f_\alpha = (f,H_\alpha)_{L^2}$ as in Subsection \ref{Ap1},
we get
\begin{align*}
\sum_{\alpha}
e^{2\tau \, \left(\log(1+\lambda_\alpha)\right)^{\frac{2}{s}}} |f_\alpha|^2
= \|e^{\tau(\log(\mathcal{H}+1))^{\frac{2}{s}}}f\|^2_{L^2(\mathbb{R}^3)}.
\end{align*}
We rephrase the previous identity as follows

\begin{align*}
\sum_{\alpha\in\mathbb{N}^3}
[h_{\tau,k}(1+\lambda_\alpha)] \, (1+\lambda_\alpha)^k |f_\alpha|^2
=\|e^{\tau(\log(\mathcal{H}+1))^{\frac{2}{s}}}f\|^2_{L^2(\mathbb{R}^3)}
\end{align*}
where
$h_{\tau,k}(x) = \frac{e^{2\tau\,\left(\log x\right)^{\frac{2}{s}}}}{x^k}.$
It is easy to check that
\begin{equation}\label{young}
\forall x\geq1, \quad h_{\tau,k}(x) \geq
e^{-\frac{2-s}{2}\big(\frac{s}{4\tau}\big)^{\frac{s}{2-s}}k^{\frac{2}{2-s}}}.
\end{equation}
Indeed, for $0<s<2$, using Young's inequality
$$xy\leq \frac{1}{p}x^p+\frac{1}{q}y^q,
\,\,\text{where}\,\,\frac{1}{p}+\frac{1}{q}=1,$$
with $p=\frac{2}{2-s}$, $q=\frac{2}{s}$,
we obtain
$$k\log x\leq\frac{2-s}{2}\Big[\Big(\frac{s}{4\tau}\Big)^{\frac{s}{2}}k\Big]^{\frac{2}{2-s}}+2\tau(\log x)^{\frac{2}{s}}.$$
Therefore,
$$h_{\tau,k}(x)=e^{2\tau\,\left(\log x\right)^{\frac{2}{s}}-k\log x}\geq e^{-\frac{2-s}{2}\Big[\Big(\frac{s}{4\tau}\Big)^{\frac{s}{2}}k\Big]^{\frac{2}{2-s}}}.$$
Then \eqref{young} follows immediately.
We conclude that
\begin{align*}
\|e^{\tau(\log(\mathcal{H}+1))^{\frac{2}{s}}}f\|^2_{L^2(\mathbb{R}^3)}
&=
\sum_{\alpha\in\mathbb{N}^3}
h_{\tau,k} (1+\lambda_{\alpha})^k \, |f_{\alpha}|^2
  \\
&\geq
e^{-\frac{2-s}{2}\big(\frac{s}{4\tau}\big)^{\frac{s}{2-s}}k^{\frac{2}{2-s}}} \,
\sum_{\alpha\in\mathbb{N}^3}
(1+\lambda_\alpha)^k \, |f_\alpha|^2=e^{-C_s\big(\frac{1}{\tau}\big)^{\frac{s}{2-s}}k^{\frac{2}{2-s}}}\|f\|^2_{Q^{k}(\mathbb{R}^3)},
\end{align*}
where we used the result of Subsection \ref{Ap1}.  This ends the proof.
\end{proof}

\bigskip
\noindent {\bf Acknowledgements.}
The second author would like to express his sincere thanks to
Prof. Chao-Jiang Xu for stimulating suggestions.
The research of the second author is supported by ``The Fundamental Research Funds for Central Universities''.


\begin{thebibliography}{99}
\bibitem{BHRV0}
J.-M. Barbaroux, D. Hundertmark, T. Ried, S. Vugalter,
Gevrey smoothing for weak solutions of the fully nonlinear homogeneous Boltzmann and Kac equations without cutoff for Maxwellian molecules, arXiv:1509.01444.

\bibitem{BHRV}
J.-M. Barbaroux, D. Hundertmark, T. Ried, S. Vugalter, Strong smoothing for the non-cutoff homogeneous Boltzmann equation for Maxwellian molecules with Debye-Yukawa type interaction, arXiv:1512.05134.


\bibitem{Boby}
A.V. Bobylev, The theory of the nonlinear spatially uniform Boltzmann equation for Maxwell
molecules, Soviet Sci. Rev. Sect. C Math. Phys. 7 (1988), 111-233.




\bibitem{Cerci}
C.~Cercignani,~~The Boltzmann Equation and its Applications,~Applied Mathematical Sciences,~vol. 67~(1988),~Springer-Verlag, New York.

\bibitem{DFT}
L.~Desvillettes,~G. Furioli,~E. Terraneo,~~Propagation of Gevrey regularity for solutions of Boltzmann
equation for Maxwellian molecules,~~Trans. Amer. Math. Soc. 361 (2009), 1731-1747.

\bibitem{DW}
L.~Desvillettes,~B.~Wennberg,~~Smoothness of the solution of the spatially homogeneous Boltzmann
equation without cutoff,~Comm.~Partial Differential Equations,~29 (2004),~no.~1-2,~133-155.



\bibitem{Dole}  
\newblock E.~Dolera,
\newblock On the computation of the spectrum of the linearized Boltzmann collision operator for Maxwellian molecules,
\newblock Boll. Unione Mat. Ital.(9), 4 (2011), 47-68.

\bibitem{G-N}
L. Glangetas, M. Najeme, Analytical regularizing effect for the radial homogeneous
     Boltzmann equation, Kinet. Relat. Models 6 (2013), no. 2, 407-427.


\bibitem{GL-2015}
L. Glangetas, H.-G. Li, Sharp regularity and Cauchy problem of the spatially homogeneous Boltzmann equation with Debye-Yukawa potential. arXiv:1512.06665.



\bibitem{GLX_2015}
L. Glangetas, H.-G. Li, C.-J. Xu, Sharp regularity properties for the non-cutoff spatially homogeneous Boltzmann equation. Kinet. Relat. Models 9 (2016), no. 2, 299-371.

\bibitem{GPR}
T. Gramchev , S. Pilipovi$\acute{c}$, L. Rodino,  Classes of degenerate elliptic operators in Gelfand-Shilov spaces. New Developments in Pseudo-Differential Operators. Birkh$\ddot{a}$user Basel, 2009: 15-31.

\bibitem{Jones}
M. N.\,Jones,  Spherical harmonics and tensors for classical field theory. UK: Research Studies Press, 1985.


\bibitem{L-X}
N. Lekrine, C.-J. Xu, Gevrey regularizing effect of the Cauchy problem for non-cutoff homogeneous Kac equation.  Kinet. Relat. Models 2 (2009), no. 4, 647-666.


\bibitem{NYKC1}
N. Lerner,~Y.~Morimoto,~K.~Pravda-Starov,~C.-J.~Xu, Spectral and phase space analysis of the linearized non-cutoff Kac collision operator. J. Math. Pures Appl.  100 (2013), no. 6, 832-867.

\bibitem{NYKC2}
N. Lerner, Y. Morimoto, K. Pravda-Starov, C.-J. Xu, Phase space analysis and functional calculus for
the linearized Landau and Boltzmann operators, Kinet. Relat. Models 6 (2013), no. 3, 625-648.

\bibitem{NYKC3}
N. Lerner,  Y. Morimoto, K. Pravda-Starov,C.-J. Xu, Gelfand-Shilov smoothing properties of the radially symmetric spatially homogeneous Boltzmann equation without angular cutoff. Journal of Differential Equations, 256 (2014), no. 2, 797-831.

\bibitem{HAOLI}
H.-G. Li,  Cauchy problem for linearized non-cutoff
Boltzmann equation with distribution initial datum. Acta Mathematica Scientia, 35B (2015), no. 2, 459-476.

\bibitem{Morimoto}
Y. Morimoto, Hypoellipticity for infinitely degenerate elliptic operators. Osaka J. Math. 24 (1987), no. 1, 13-35.

\bibitem{Morimoto2}
Y. Morimoto, A remark on Cannone-Karch solutions to the homogeneous Boltzmann equation
for Maxwellian molecules, Kinetic and Related Models, 5 (2012), no. 3, 551-561.

\bibitem{MPX}
Y. Morimoto, K. Pravda-Starov, C.-J. Xu,  A remark on the ultra-analytic smoothing properties of the spatially homogeneous Landau equation. Kinetic and Related Models 6 (2013), no. 4, 715-727.

\bibitem{MU}
Y.~Morimoto,~S. Ukai,~~Gevrey smoothing effect of solutions for spatially homogeneous nonlinear Boltzmann equation without angular cutoff,~J. Pseudo-Differ. Oper., Appl. 1 (2010), no. 1, 139-159.

\bibitem{YSCT}
Y. Morimoto,~S. Ukai,~C.-J. Xu,~T. Yang, Regularity of solutions to the spatially homogeneous Boltzmann equation without angular cutoff. Discrete Contin. Dyn. Syst. 24 (2009), no. 1, 187-212.

\bibitem{Morimoto-Xu1}
Y. Morimoto,  C.-J. Xu,
Logarithmic Sobolev inequality and semi-linear Dirichlet problems for infinitely degenerate elliptic operators.
Ast\'erisque, No. 284 (2003), 245-264.

\bibitem{Morimoto-Xu2}
Y. Morimoto,  C.-J. Xu,
Nonlinear hypoellipticity of infinite type.
Funkcial. Ekvac. 50 (2007), no. 1, 33-65.



\bibitem{JCSlater}
J. C. Slater,  Quantum theory of atomic structure. Vol. 1. New York: McGraw-Hill, 1960.

\bibitem{San}
G.~Sansone,~Orthogonal Functions.~Pure and Applied Mathematics.~Vol.~IX.~1959~Interscience Publishers,~New York.~~Reprinted by Dover~Publications~1991.

\bibitem{Shubin}
M. Shubin, Pseudodifferential Operators and Spectral theory. Springer Series in Soviet
Mathematics, Springer-Verlag, Berlin, 1987.

\bibitem{Ukai}
S. Ukai,~~Local solutions in Gevrey classes to the nonlinear Boltzmann equation without cutoff,~Japan.J. Appl. Math. 1 (1984), no.1,~141-156.

\bibitem{Villani}
C. Villani,  A review of mathematical topics in collisional kinetic theory. Handbook of mathematical fluid dynamics, 2002, 1: 71-305.

\bibitem{WU}
C.S.~Wang Chang,~G.E.~Uhlenbeck,~~On the propagation of sound in monoatomic gases,~~Univ.
of Michigan Press.~Ann Arbor, Michigan.~~Reprinted in 1970 in Studies in Statistical Mechanics.
Vol.~V.~Edited by J.L.~Lebowitz and E.~Montroll,~North-Holland


\bibitem{TZ}
T.-F. Zhang, Z. Yin,~~Gevrey regularity of spatially homogeneous Boltzmann equation without cutoff,
J. Differential Equations,~ 253 (2012), no. 4, 1172-1190.

\end{thebibliography}
\end{document}